\documentclass[12pt]{article}
\usepackage{amsmath,amssymb,amsthm,amscd,graphicx}
\usepackage{latexsym}
\usepackage{enumerate}
\usepackage{amsbsy, amsfonts, textcomp}
\newtheorem{theorem}{Theorem}[]
\newtheorem{proposition}[theorem]{Proposition}

\newtheorem{conjecture}[theorem]{Conjecture}
\theoremstyle{definition}

\newtheorem{example}[theorem]{Example}

\theoremstyle{remark}
\newtheorem{remark}[theorem]{Remark}

\newcommand{\C}{\mathbb{C}}

\def \max{\operatorname{max}}

\title{Jacobian conjecture and nilpotency}

\author{El\.zbieta Adamus\footnote{Faculty of Applied Mathematics,  AGH University of Science and
Technology, al. Mickiewicza 30, 30-059 Krak\'ow, Poland, e-mail:
esowa@agh.edu.pl}, Pawe\l \ Bogdan\footnote{Faculty of Mathematics
and Computer Science, Jagiellonian University, ul. \L ojasiewicza
6, 30-348 Krak\'ow, Poland, e-mail: pawel.bogdan@uj.edu.pl},
Teresa Crespo\footnote{Departament d'\`{A}lgebra i Geometria,
Universitat de Barcelona, Gran Via de les Corts Catalanes 585,
08007 Barcelona, Spain,  e-mail: teresa.crespo@ub.edu} and
Zbigniew Hajto\footnote{Faculty of Mathematics and Computer
Science, Jagiellonian University, ul. \L ojasiewicza 6, 30-348
Krak\'ow, Poland, e-mail: zbigniew.hajto@uj.edu.pl}}

\date{\today}

\begin{document}
\maketitle

\begin{abstract}
For $K$ a field of characteristic 0 and $d$ any integer
number $\geq 2$, we prove the invertibility of polynomial maps $F:K^d
\rightarrow K^d$ of the form $F=Id+H$, where each $H_i$ is the
cube of a linear form and the Jacobian matrix $JH$ satisfies
$(JH)^3=0$. Our proof uses the inversion algorithm for polynomial
maps presented in our previous paper. Our current result leads us
to formulate a conjecture relating the nilpotency degree of the
matrix $JH$ with the number of necessary steps in the inversion
algorithm.

\noindent {\bf Keywords:} Polynomial automorphism, Jacobian
Problem, inversion algorithm.

\noindent {\bf MSC2010:} 14R15, 14R10.
\end{abstract}

\section{Introduction}

The Jacobian Conjecture originated in the question raised by Keller in \cite{K} on the invertibility of polynomial maps with Jacobian determinant equal to~1. The question is still open in spite of the efforts of many mathematicians. We recall in the sequel the precise statement of the Jacobian Conjecture, some reduction theorems and other results we shall use. We refer to \cite{E} for a detailed account of the research on the Jacobian Conjecture and related topics.

Let $K$ be a field and $K[X]=K[X_1,\dots,X_d]$ the polynomial ring in the variables $X_1,\dots,X_d$ over $K$. A \emph{polynomial map} is a map $F=(F_1,\dots,F_d):K^d \rightarrow K^d$ of the form

$$(X_1,\dots,X_d)\mapsto (F_1(X_1,\dots,X_d),\dots,F_d(X_1,\dots,X_d)),$$

\noindent where $F_i \in K[X], 1 \leq i \leq d$. The polynomial map $F$ is \emph{invertible} if there exists a polynomial map $G=(G_1,\dots,G_d):K^d \rightarrow K^d$ such that $X_i=G_i(F_1,\dots,F_d),  1 \leq i \leq d$. We shall call $F$ a \emph{Keller map} if the Jacobian matrix

$$JF=\left(\dfrac{\partial F_i}{\partial X_j}\right)_{\substack{1\leq i \leq d \\ 1\leq j \leq d}}$$

\noindent has determinant equal to 1. Clearly an invertible polynomial map $F$ has a Jacobian matrix $JF$ with non zero determinant and may be transformed into a Keller map by composition with the linear automorphism with matrix $JF(0)^{-1}$.

\vspace{0.5cm}
\noindent {\bf Jacobian Conjecture.} {\it Let $K$ be a field of characteristic zero. A Keller map $F:K^d \rightarrow K^d$ is invertible.}

\vspace{0.5cm}
In the sequel, $K$ will always denote a {\bf field of characteristic zero.}

For $F=(F_1,\dots,F_d) \in K[X]^d$, we define the \emph{degree} of $F$ as $\deg F= \max \{\deg F_i : 1\leq i \leq d\}$. It is known that if $F$ is a polynomial automorphism, then $\deg F^{-1} \leq (\deg F)^{d-1}$ (see \cite{BCW} or \cite{RW}).

The Jacobian conjecture for quadratic maps was proved by Wang in \cite{W}.  We state now the reduction of the Jacobian conjecture to the case of maps of third degree (see \cite{BCW}, \cite{Y}, \cite{D1} and \cite{D2}).

\begin{proposition}\label{red} \begin{enumerate}[a)] \item (Bass-Connell-Wright-Yagzhev) Given a Keller map $F:\C^d \rightarrow \C^d$, there exists a Keller map $\widetilde{F}:\C^D \rightarrow \C^D$, $D\geq d$ of the form $\widetilde{F}=Id+H$, where $H(X)$ is a homogeneous cubic map and having the following property: if $\widetilde{F}$ is invertible, then $F$ is invertible too.
\item (Dru\.zkowski) The cubic part $H$ may be chosen of the form

$$\left( (\sum_{j=1}^D a_{1j} X_j)^3, \dots,(\sum_{j=1}^D a_{Dj} X_j)^3 \right)$$

\noindent and with the matrix $A=(a_{ij})_{\substack{1\leq i \leq D \\ 1\leq j \leq D}}$ satisfying $A^2=0$.

\end{enumerate}

\end{proposition}

Let us note that for a polynomial map in $d$ variables of the form $F=Id+H$, with $H$ homogeneous of degree at least two, the condition $det(JF)=1$ is equivalent to $(JH)^d=0$. Given a polynomial map $F:\C^d \rightarrow \C^d$, we shall call a polynomial $P \in \C[X_1,\dots,X_d]$ \emph{invariant under $F$} if $P(F_1,\dots,F_d)=P(X_1\dots,X_d)$. A polynomial map $F:\C^d \rightarrow \C^d$ of the form $F=Id+H$ is called a \emph{quasi-translation} if $F^{-1}=Id-H$.

In \cite{Alg} we proved the following algorithmic equivalent statement to the Jacobian conjecture for homogeneous polynomial maps. Although it is there stated in the complex case, it is clearly valid for any field $K$ of characteristic 0.

\begin{theorem}\label{theo} Let $F:K^d \rightarrow K^d$ be a polynomial map of the form

$$\left\{ \begin{array}{lll} F_1(X_1,\dots,X_d)&=& X_1+H_1(X_1,\dots,X_d) \\ F_2(X_1,\dots,X_d)&=& X_2+H_2(X_1,\dots,X_d)\\ & \vdots & \\  F_d(X_1,\dots,X_d)&=& X_d+H_d(X_1,\dots,X_d) , \end{array} \right.$$

\noindent where $H_i(X_1,\dots,X_d)$ is a homogeneous polynomial in $X_1,\dots,X_d$, of degree~$3$, $1\leq i \leq d$.  For each $i=1,\dots,d$, we consider the polynomial sequence $(P_j^i)$ defined in the following way

$$\begin{array}{lll} P_0^i(X_1,\dots,X_d) &= & X_i, \\
P_1^i(X_1,\dots,X_d) &=& H_i , \end{array}$$

\noindent and, assuming $P_{j-1}^i$ is defined,

$$P_j^i(X_1,\dots,X_d) = P_{j-1}^i(F_1,\dots,F_d)-P_{j-1}^i(X_1,\dots,X_d).$$

Then $F$ is invertible if and only if for all $i=1,\dots, d$, there exists an integer $m_i$ such that $P_{m_i}^i = 0$. In this case, the inverse map $G$ of $F$ is given by

$$G_i(Y_1,Y_2,\dots,Y_d)= \sum_{l=0}^{m_i-1} (-1)^l P_l^i(Y_1,Y_2,\dots,Y_d), \, 1 \leq i \leq d.$$

\end{theorem}

In this paper we shall use this theorem to prove that a polynomial map of the Dru\.zkowski form such that the Jacobian matrix of $H$ is nilpotent of degree at most 3 is invertible.

\section{Main result}

We shall prove the following theorem

\begin{theorem}\label{theo2}
Let $F:K^d \rightarrow K^d$ be a polynomial map of the form

$$\left\{ \begin{array}{lll} F_1(X_1,\dots,X_d)&=& X_1+H_1(X_1,\dots,X_d) \\ F_2(X_1,\dots,X_d)&=& X_2+H_2(X_1,\dots,X_d)\\ & \vdots & \\  F_d(X_1,\dots,X_d)&=& X_d+H_d(X_1,\dots,X_d) , \end{array} \right.$$

\noindent where $H_i(X_1,\dots,X_d)=L_i(X_1,\dots,X_d)^3$ and $L_i(X_1,\dots,X_d)=a_{i1}X_1+\dots +a_{id} X_d$, $1\leq i \leq d$. We consider the jacobian matrix $JH$ of $H:=(H_1,\dots,H_d)$

$$JH= \left( \begin{array}{ccc} \dfrac{\partial H_1}{\partial X_1} & \dots &  \dfrac{\partial H_1}{\partial X_d}\\ & \vdots & \\ \dfrac{\partial H_d}{\partial X_1} & \dots &  \dfrac{\partial H_d}{\partial X_d} \end{array} \right). $$

\noindent If $(JH)^3=0$, then $F$ is invertible and $\deg F^{-1} \leq 9$.
\end{theorem}

The theorem follows from Theorem \ref{theo} and the Proposition
\ref{prop} below.

\begin{proposition}\label{prop} Let $F$ be as in Theorem \ref{theo2} and let the polynomials $P_j^i$ be defined as in Theorem \ref{theo}. If $(JH)^3=0$, then  $\forall i=1,\dots,d, \, \deg P_j^i \leq 9, \, \text{for} \, j=2,3,4, \, \text{and} \,P_5^i=0$.
\end{proposition}

\noindent {\it Proof of Theorem \ref{theo2} assuming Proposition \ref{prop}.} If $(JH)^3=0$, we have $P_5^i=0$, for all $i$, by Proposition \ref{prop} and this implies, by Theorem \ref{theo}, that $F$ is invertible and the inverse map $G$ of $F$ is given by $G_i(Y_1,Y_2,\dots,Y_d)= Y_i-P_1^i(Y_1,Y_2,\dots,Y_d)+P_2^i(Y_1,Y_2,\dots,Y_d)-P_3^i(Y_1,Y_2,\dots,Y_d)+P_4^i(Y_1,Y_2,\dots,Y_d)$. By the bound on the degrees given by Proposition \ref{prop}, we obtain $\deg F^{-1} \leq 9$.

\vspace{0.5cm}
\noindent {\it Proof of Proposition \ref{prop}.}
This proof is purely computational. We shall use the expression of the polynomials $P_j^i$ as a sum of homogeneous polynomials obtained in \cite{Alg}, proof of Theorem 4. From the condition $(JH)^3=0$, we shall derive some equations which will allow to prove the vanishing of some homogeneous summands of the polynomials $P_j^i$.

Since  $(JH)^3=0$, we have

$$
\sum_{j,k=1}^d \dfrac{\partial H_i}{\partial X_j}\dfrac{\partial H_j}{\partial X_k}\dfrac{\partial H_k}{\partial X_l}=0, \quad \forall \, i, l =1,\dots,d.
$$

\noindent i.e.

$$
\sum_{j,k=1}^d (3L_i^2 a_{ij}) (3L_j^2 a_{jk}) (3L_k^2 a_{kl}) =27 L_i^2 (\sum_{j,k=1}^d  a_{ij}L_j^2 a_{jk} L_k^2 a_{kl})=0, \quad \forall \, i, l =1,\dots,d.
$$

\noindent which implies

\begin{equation}\label{eq3}
\sum_{j,k=1}^d a_{ij}  a_{jk} a_{kl} L_j^2 L_k^2  =0, \quad \forall \, i, l =1,\dots,d,
\end{equation}

\noindent since $L_i=0$ implies $a_{ij}=0$, for all $j$, hence (\ref{eq3}). Applying $\dfrac{\partial}{\partial X_m}$, we obtain

\begin{equation}\label{eq4}
\sum_{j,k=1}^d a_{ij} a_{jk} a_{kl}(a_{jm} L_jL_k^2+ a_{km}L_j^2L_k ) =0, \quad \forall \, i, l, m =1,\dots,d.
\end{equation}

\noindent Further, from (\ref{eq3}), we obtain

\begin{equation}\label{eq8}
\sum_{j,k=1}^d a_{ij} a_{jk}L_j^2L_k^3 =0, \quad \forall \, i=1,\dots,d
\end{equation}

\noindent since $\sum_{j,k=1}^d a_{ij} a_{jk}L_j^2L_k^3= \sum_{j,k=1}^d a_{ij} a_{jk}L_j^2L_k^2 (\sum_{l=1}^d a_{kl}X_l)= \sum_{l=1}^d (\sum_{j,k=1}^d a_{ij}  a_{jk} a_{kl} L_j^2 L_k^2) X_l$. Now, from (\ref{eq4}), we obtain similarly

$$\sum_{j,k=1}^d a_{ij} a_{jk} (a_{jm}L_jL_k^3+a_{km}L_j^2L_k^2 ) =0, \quad \forall \, i, m =1,\dots,d$$

\noindent and, using (\ref{eq3})

\begin{equation}\label{eq9}
\sum_{j,k=1}^d a_{ij} a_{jk} a_{jm} L_jL_k^3 =0, \quad \forall \, i, m =1,\dots,d
\end{equation}

\noindent Applying $\partial/\partial X_n$ to (\ref{eq9}), we obtain

$$\sum_{j,k=1}^d a_{ij} a_{jk} a_{jm} (a_{jn} L_k^3+3a_{kn} L_j L_k^2) =0$$

\noindent i.e.

\begin{equation}\label{eq10}
\sum_{j,k=1}^d a_{ij} a_{jk} a_{jm} a_{jn} L_k^3= - 3\sum_{j,k=1}^d a_{ij} a_{jk} a_{jm} a_{kn} L_j L_k^2
\end{equation}

We shall use equations (\ref{eq8}), (\ref{eq9}) and (\ref{eq10}) repeatedly to prove the vanishing of some homogeneous summands of the polynomials $P_2^i, P_3^i, P_4^i$ and the vanishing of $P_5^i$. Since the calculations are very similar we will detail them only the first time we apply each of these equations.

\noindent We have

$$P_2^i=H_i(F)-H_i(X)=H_i(X+H)-H_i(X)= Q_{21}^i+Q_{22}^i+Q_{23}^i
$$

\noindent where

$$\begin{array}{lll} Q_{21}^i &=& \sum_{j=1}^d H_j \dfrac{\partial H_i}{\partial X_j}= 3L_i^2 \sum_{j=1}^d a_{ij} L_j^3  \\ [10pt] Q_{22}^i &=& \dfrac 1 2 \sum_{1\leq j,k \leq d} H_j H_k \dfrac{\partial^2 H_i}{\partial X_j\partial X_k}=3L_i \sum_{1\leq j,k \leq d} a_{ij} a_{ik} L_j^3 L_k^3 \\ [10pt] Q_{23}^i &=& \dfrac 1 6 \sum_{1\leq j,k,l \leq d} H_jH_kH_l \dfrac{\partial^3 H_i}{\partial X_j \partial X_k \partial X_l} = \sum_{1\leq j,k,l \leq d}  a_{ij} a_{ik} a_{il} L_j^3 L_k^3 L_l^3, \end{array} $$

\noindent hence, in particular $\deg P_2^i \leq 9$.

In order to determine $P_3^i$, we need to compute the derivatives
of $Q_{21}^i, Q_{22}^i$ and $Q_{23}^i$. We compute first the
derivatives of $Q_{21}^i$.

$$\begin{array}{rcl} \dfrac{\partial Q_{21}^i}{\partial X_k} &=& 6 a_{ik}L_i \sum_{j=1}^d  a_{ij}L_j^3+9L_i^2 \sum_{j=1}^d  a_{ij}a_{jk}L_j^2\\ [12pt]
\dfrac{\partial^2 Q_{21}^i}{\partial X_k\partial X_l} &=& 6 a_{ik} a_{il} \sum_{j=1}^d  a_{ij}L_j^3+18L_i a_{ik} \sum_{j=1}^d  a_{ij}a_{jl}L_j^2 \\ && +18L_i a_{il} \sum_{j=1}^d  a_{ij}a_{jk}L_j^2+18L_i^2 \sum_{j=1}^d  a_{ij}a_{jk}a_{jl}L_j \end{array}$$

\noindent The derivative $\partial^3 Q_{21}^i /\partial X_k\partial X_l\partial X_m$ is a sum of terms of one  of the following forms

$$ 18 a_{ik} a_{il} \sum_{j=1}^d  a_{ij}a_{jm}L_j^2, \quad 36L_i a_{ik} \sum_{j=1}^d  a_{ij}a_{jl}a_{jm}L_j,$$

\noindent up to a permutation of the set $\{k,l,m\}$, and the summand $18L_i^2 \sum_{j=1}^d  a_{ij}a_{jk}a_{jl}a_{jm}$. The derivative $\partial^4 Q_{21}^i /\partial X_k\partial X_l\partial X_m\partial X_n $ is a sum of terms of one of the following forms

$$36 a_{ik} a_{il} \sum_{j=1}^d a_{ij}a_{jm}a_{jn} L_j, \quad  36L_i a_{ik} \sum_{j=1}^d a_{ij}a_{jl}a_{jm}a_{jn},$$

\noindent up to a permutation of the set $\{k,l,m,n\}$. The derivative $\partial^5 Q_{21}^i /\partial X_k\partial X_l\partial X_m \partial X_n \partial X_p$ is a sum of terms of the form

$$ 36 a_{ik} a_{il} \sum_{j=1}^d  a_{ij}a_{jm}a_{jn}a_{jp}$$

\noindent up to a permutation of the set $\{k,l,m,n,p\}$.

We compute now the derivatives of $Q_{22}^i$.

$$\begin{array}{rcl} \dfrac{\partial Q_{22}^i}{\partial X_l} &=& 3a_{il} \sum_{1\leq j,k \leq d}  a_{ij} a_{ik}L_j^3 L_k^3+9L_i \sum_{1\leq j,k \leq d}  a_{ij} a_{ik}a_{jl}L_j^2 L_k^3\\ && +9L_i \sum_{1\leq j,k \leq d}  a_{ij} a_{ik}a_{kl} L_j^3 L_k^2\end{array}$$

\noindent The derivative $\partial^2 Q_{22}^i/ \partial X_l\partial X_m$ is a sum of terms of one of the following forms

$$\begin{array}{ll} 9a_{il} \sum_{1\leq j,k \leq d}  a_{ij} a_{ik}a_{jm}L_j^2 L_k^3, & 18L_i \sum_{1\leq j,k \leq d} a_{ij} a_{ik}a_{jl}a_{jm} L_j L_k^3 \\  27 L_i \sum_{1\leq j,k \leq d}  a_{ij} a_{ik}a_{jl}a_{km} L_j^2 L_k^2, & \end{array}$$

\noindent up to switching $l$ with $m$, and $j$ with $k$. The derivative $\partial^3 Q_{22}^i / \partial X_l\partial X_m\partial X_n$ is a sum of terms of one of the following forms

$$\begin{array}{ll} 18a_{il} \sum_{1\leq j,k \leq d}  a_{ij} a_{ik}a_{jm}a_{jn} L_j L_k^3, & 27a_{il} \sum_{1\leq j,k \leq d}  a_{ij} a_{ik}a_{jm} a_{kn} L_j^2 L_k^2,\\ 18L_i \sum_{1\leq j,k \leq d} a_{jn}  a_{ij} a_{ik}a_{jl}a_{jm}L_k^3, & 54 L_i \sum_{1\leq j,k \leq d} a_{ij} a_{ik} a_{jl} a_{jm} a_{kn} L_j L_k^2, \end{array}$$

\noindent up to a permutation of the set $\{l,m,n\}$ and a switch of $j$ with $k$. The derivative $\partial^4 Q_{22}^i/ \partial X_l\partial X_m\partial X_n\partial X_p$ is a sum of terms of one of the following forms

$$\begin{array}{ll} 18a_{il} \sum_{1\leq j,k \leq d}   a_{ij} a_{ik}a_{jm}a_{jn}a_{jp} L_k^3,& 54 a_{il} \sum_{1\leq j,k \leq d}  a_{ij} a_{ik}a_{jm}a_{jn}a_{kp} L_j L_k^2, \\ 54L_i \sum_{1\leq j,k \leq d}   a_{ij} a_{ik}a_{jl}a_{jm}a_{jn}a_{kp} L_k^2,& 108L_i \sum_{1\leq j,k \leq d}  a_{ij} a_{ik}a_{jl}a_{jm}a_{kn}a_{kp}L_j L_k, \end{array}$$

\noindent up to a permutation of the set $\{l,m,n,p\}$ and a switch of $j$ with $k$. The derivative $\partial^5 Q_{22}^i/\partial X_l\partial X_m\partial X_n\partial X_p\partial X_q$ is a sum of terms of one of the following forms

$$\begin{array}{ll}  54 a_{il} \sum_{1\leq j,k \leq d}  a_{ij} a_{ik}a_{jm}a_{jn}a_{jp}a_{kq}L_k^2, &  108 a_{il} \sum_{1\leq j,k \leq d}  a_{ij} a_{ik}a_{jm}a_{jn}a_{kp}a_{kq} L_j L_k, \\ 108 L_i \sum_{1\leq j,k \leq d}   a_{ij} a_{ik}a_{jl}a_{jm}a_{jn}a_{kp}a_{kq}L_k, &
\end{array}
$$

\noindent up to a permutation of the set $\{l,m,n,p,q\}$ and a switch of $j$ with $k$. The derivative $\partial^6 Q_{22}^i/\partial X_l\partial X_m\partial X_n\partial X_p\partial X_q \partial X_r$ is a sum of terms of one of the following forms

$$  108 a_{il} \sum_{1\leq j,k \leq d}  a_{ij} a_{ik}a_{jm}a_{jn}a_{jp}a_{kq} a_{kr} L_k, \quad 108 L_i \sum_{1\leq j,k \leq d}   a_{ij} a_{ik}a_{jl}a_{jm}a_{jn}a_{kp}a_{kq}a_{kr},
$$

\noindent up to a permutation of the set $\{l,m,n,p,q,r\}$ and a switch of $j$ with $k$. The derivative $\partial^7 Q_{22}^i/\partial X_l\partial X_m\partial X_n\partial X_p\partial X_q \partial X_r\partial X_s$ is a sum of terms of the form

$$  108 a_{il} \sum_{1\leq j,k \leq d}  a_{ij} a_{ik}a_{jm}a_{jn}a_{jp}a_{kq} a_{kr} a_{ks},
$$

\noindent up to a permutation of the set $\{l,m,n,p,q,r,s\}$ and a switch of $j$ with $k$.

Finally we compute the derivatives of $Q_{23}^i$.

$$\dfrac{\partial Q_{23}^i}{\partial X_m}  = 3\sum_{j,k,l}  a_{ij} a_{ik} a_{il}a_{jm}L_j^2 L_k^3 L_l^3+3\sum_{j,k,l}  a_{ij} a_{ik} a_{il}a_{km}L_j^3 L_k^2 L_l^3+3\sum_{j,k,l}  a_{ij} a_{ik} a_{il}a_{lm}L_j^3 L_k^3 L_l^2.$$

\noindent The derivative $\partial^2 Q_{23}^i / \partial X_m\partial X_n$ is a sum of terms of one of the following forms

$$ 6\sum_{j,k,l}  a_{ij} a_{ik} a_{il}a_{jm}a_{jn}L_j L_k^3 L_l^3, \quad 9\sum_{j,k,l}  a_{ij} a_{ik} a_{il}a_{jm}a_{kn}L_j^2 L_k^2 L_l^3$$

\noindent up to a permutation of the set $\{m,n\}$ and a permutation of the set $\{i,j,k\}$. The derivative $\partial^3 Q_{23}^i / \partial X_m\partial X_n\partial X_p$   is a sum of terms of one of the following forms

$$\begin{array}{ll} 6\sum_{j,k,l}   a_{ij} a_{ik} a_{il}a_{jm}a_{jn}a_{jp} L_k^3 L_l^3, &  18\sum_{j,k,l}  a_{ij} a_{ik} a_{il}a_{jm}a_{jn}a_{kp}L_j L_k^2 L_l^3\\
                        27 \sum_{j,k,l}  a_{ij} a_{ik} a_{il}a_{jm}a_{kn}a_{lp}L_j^2 L_k^2 L_l^2, & \end{array}$$

\noindent up to a permutation of the set $\{m,n,p\}$ and a permutation of the set $\{i,j,k\}$. The derivative $\partial^4 Q_{23}^i / \partial X_m\partial X_n\partial X_p\partial X_q$   is a sum of terms of one of the following forms

$$\begin{array}{ll} 18\sum_{j,k,l}   a_{ij} a_{ik} a_{il}a_{jm}a_{jn}a_{jp} a_{kq} L_k^2 L_l^3, & 36 \sum_{j,k,l}  a_{ij} a_{ik} a_{il}a_{jm}a_{jn}a_{kp}a_{kq} L_j L_k L_l^3\\54\sum_{j,k,l}  a_{ij} a_{ik} a_{il}a_{jm}a_{jn}a_{kp}a_{lq} L_j L_k^2 L_l^2, & \end{array}
$$

\noindent up to a permutation of the set $\{m,n,p,q\}$ and a permutation of the set $\{i,j,k\}$.  The derivative $\partial^5 Q_{23}^i / \partial X_m\partial X_n\partial X_p\partial X_q\partial X_r$   is a sum of terms of one of the following forms

$$\begin{array}{ll} 36\sum_{j,k,l}   a_{ij} a_{ik} a_{il}a_{jm}a_{jn}a_{jp} a_{kq} a_{kr} L_k L_l^3, & 54 \sum_{j,k,l}  a_{ij} a_{ik} a_{il}a_{jm}a_{jn}a_{jp}a_{kq} a_{lr} L_k^2 L_l^2\\108 \sum_{j,k,l}  a_{ij} a_{ik} a_{il}a_{jm}a_{jn}a_{kp}a_{kq} a_{lr} L_j L_k L_l^2, & \end{array}
$$

\noindent up to a permutation of the set $\{m,n,p,q,r\}$ and a permutation of the set $\{i,j,k\}$.  The derivative $\partial^6 Q_{23}^i / \partial X_m\partial X_n\partial X_p\partial X_q\partial X_r\partial X_s$   is a sum of terms of one of the following forms

$$\begin{array}{ll} 36\sum_{j,k,l}   a_{ij} a_{ik} a_{il}a_{jm}a_{jn}a_{jp} a_{kq} a_{kr} a_{ks} L_l^3, & 108 \sum_{j,k,l}  a_{ij} a_{ik} a_{il}a_{jm}a_{jn}a_{jp}a_{kq} a_{lr} a_{ks} L_k L_l^2\\216 \sum_{j,k,l}  a_{ij} a_{ik} a_{il}a_{jm}a_{jn}a_{kp}a_{kq} a_{lr} a_{ls} L_j L_k L_l, & \end{array}
$$

\noindent up to a permutation of the set $\{m,n,p,q,r,s\}$ and a permutation of the set $\{i,j,k\}$. The derivative $\partial^7 Q_{23}^i / \partial X_m\partial X_n\partial X_p\partial X_q\partial X_r\partial X_s\partial X_t$   is a sum of terms of one of the following forms

$$108 \sum_{j,k,l}   a_{ij} a_{ik} a_{il}a_{jm}a_{jn}a_{jp} a_{kq} a_{kr} a_{ks}a_{lt} L_l^2, \quad 216 \sum_{j,k,l}  a_{ij} a_{ik} a_{il}a_{jm}a_{jn}a_{jp}a_{kq} a_{lr} a_{ks} a_{lt} L_k L_l
$$

\noindent up to a permutation of the set $\{m,n,p,q,r,s,t\}$ and a permutation of the set $\{i,j,k\}$. The derivative $\partial^8 Q_{23}^i / \partial X_m\partial X_n\partial X_p\partial X_q\partial X_r\partial X_s\partial X_t \partial X_u$   is a sum of terms of the following form

$$216 \sum_{j,k,l}   a_{ij} a_{ik} a_{il}a_{jm}a_{jn}a_{jp} a_{kq} a_{kr} a_{ks}a_{lt} a_{lu} L_l
$$

\noindent up to a permutation of the set $\{m,n,p,q,r,s,t,u\}$ and a permutation of the set $\{i,j,k\}$. The derivative $\partial^9 Q_{23}^i / \partial X_m\partial X_n\partial X_p\partial X_q\partial X_r\partial X_s\partial X_t \partial X_u\partial X_v$   is a sum of terms of the following form

$$216 \sum_{j,k,l}   a_{ij} a_{ik} a_{il}a_{jm}a_{jn}a_{jp} a_{kq} a_{kr} a_{ks}a_{lt} a_{lu} a_{lv}
$$

\noindent up to a permutation of the set $\{m,n,p,q,r,s,t,u,v\}$ and a permutation of the set $\{i,j,k\}$.

 We write $P_3^i= \sum_{k=1}^{11} Q_{3k}^i$, with  $Q_{3k}^i$ homogeneous of degree $5+2k$ (see \cite{Alg} proof of Theorem 4). We have, using (\ref{eq8}),

$$\begin{array}{ll} Q_{31}^i = \sum_{k=1}^d \dfrac{\partial Q_{21}^i}{\partial X_k}H_k & = 6L_i \sum_{j,k=1}^d a_{ik}   a_{ij}L_j^3L_k^3+9L_i^2 \sum_{j,k=1}^d  a_{ij}a_{jk} L_j^2 L_k^3 \\ & = 6L_i \sum_{j,k=1}^d a_{ij}a_{ik} L_j^3 L_k^3.\end{array}$$

\noindent Using  (\ref{eq8}) and (\ref{eq9}), we obtain

$$Q_{32}^i= \dfrac 1 2 \sum_{k,l}  \dfrac{\partial^2 Q_{21}^i}{\partial X_k \partial X_l} H_kH_l + \sum_{l} \dfrac{\partial Q_{22}^i}{\partial X_l} H_l = 6 \sum_{j,k,l} a_{ij} a_{ik} a_{il} L_j^3L_k^3L_l^3,$$

\noindent since

$$\begin{array}{l} \sum_{j,k,l} a_{ij} a_{ik} a_{jl} L_j^2 L_k^3 L_l^3 =\sum_k a_{ik} L_k^3 (\sum_{j,l} a_{ij} a_{jl} L_j^2 L_l^3)\stackrel{(\ref{eq8})}{=} 0 \\ \sum_{j,k,l} a_{ij} a_{jk} a_{il} L_j^2 L_k^3 L_l^3 =\sum_l a_{il} L_l^3 (\sum_{j,k} a_{ij} a_{jk} L_j^2 L_k^3)\stackrel{(\ref{eq8})}{=} 0 \\\sum_{j,k,l} a_{ij} a_{jk} a_{jl} L_j L_k^3 L_l^3 =\sum_l L_l^3 (\sum_{j,k} a_{ij} a_{jk} a_{jl} L_j L_k^3)\stackrel{(\ref{eq9})}{=} 0 \\ \sum_{j,k,l} a_{ij} a_{ik} a_{kl} L_j^3 L_k^2 L_l^3 =\sum_j a_{ij} L_j^3 (\sum_{k,l} a_{ik} a_{kl} L_k^2 L_l^3)\stackrel{(\ref{eq8})}{=} 0 \end{array} $$

\noindent Using again (\ref{eq8}) and (\ref{eq9}), we obtain

$$\begin{array}{lll} Q_{33}^i&=& \dfrac 1 6 \sum_{k,l,m} \dfrac{\partial^3 Q_{21}^i}{\partial X_k \partial X_l \partial X_m} H_kH_lH_m+\dfrac 1 2 \sum_{l,m} \dfrac{\partial^2 Q_{22}^i}{ \partial X_l \partial X_m} H_lH_m+\sum_{m} \dfrac{\partial Q_{23}^i}{\partial X_m} H_m \\ &=& 0.
\end{array}$$

\noindent We have now, using (\ref{eq8}), (\ref{eq9}) and (\ref{eq10}),

$$\begin{array}{l} Q_{34}^i=\dfrac 1 {4!} \sum_{k,l,m,n} \dfrac{\partial^4 Q_{21}^i}{\partial X_k \partial X_l \partial X_m\partial X_n} H_kH_lH_mH_n+\dfrac 1 {3!} \sum_{l,m,n} \dfrac{\partial^3 Q_{22}^i}{ \partial X_l \partial X_m\partial X_n} H_lH_mH_n\\ +\dfrac 1 2 \sum_{m,n} \dfrac{\partial^2 Q_{23}^i}{\partial X_m\partial X_n} H_m H_n =0. \end{array}$$

\noindent We detail the use of (\ref{eq10}):

$$\begin{array}{lll} \sum_{j,k,l,m,n} a_{ij} a_{ik} a_{jl} a_{jm} a_{jn} L_k^3 L_l^3 L_m^3 L_n^3 &=& \sum_{k,m,n} a_{ik}L_k^3 L_m^3 L_n^3(\sum_{j,l} a_{ij} a_{jl} a_{jm} a_{jn} L_l^3)\\ &\stackrel{(\ref{eq10})}{=} &-3 \sum_{k,m,n} a_{ik}L_k^3 L_m^3 L_n^3(\sum_{j,l} a_{ij} a_{jl} a_{jm} a_{ln} L_j L_l^2)\\ &=& -3 \sum_{k,l,n} a_{ik}a_{jl} L_k^3 L_l^2 L_n^3(\sum_{j,m} a_{ij} a_{jl} a_{jm}  L_j L_m^3) \stackrel{(\ref{eq9})}{=}0\end{array}$$

\noindent Similarly, using again (\ref{eq8}), (\ref{eq9}) and (\ref{eq10}),

$$\begin{array}{lll} Q_{35}^i &=& \dfrac 1 {5!} \sum_{k,l,m,n,p} \dfrac{\partial^5 Q_{21}^i}{\partial X_k \partial X_l \partial X_m\partial X_n\partial X_p} H_kH_lH_mH_n H_p\\ && +\dfrac 1 {4!} \sum_{l,m,n,p} \dfrac{\partial^4 Q_{22}^i}{ \partial X_l \partial X_m\partial X_n\partial X_p} H_lH_mH_nH_p \\ && +\dfrac 1 {3!} \sum_{m,n,p} \dfrac{\partial^3 Q_{23}^i}{\partial X_m\partial X_n\partial X_p} H_m H_nH_p =0
\end{array}$$

$$\begin{array}{l} Q_{36}^i=\dfrac 1 {5!} \sum_{l,m,n,p,q} \dfrac{\partial^5 Q_{22}^i}{ \partial X_l \partial X_m\partial X_n\partial X_p\partial X_q} H_lH_mH_nH_pH_q \\ +\dfrac 1 {4!} \sum_{m,n,p,q} \dfrac{\partial^4 Q_{23}^i}{\partial X_m\partial X_n\partial X_p\partial X_q} H_m H_nH_pH_q =0
\end{array}$$

$$\begin{array}{l} Q_{37}^i=\dfrac 1 {6!} \sum_{l,m,n,p,q,r} \dfrac{\partial^6 Q_{22}^i}{ \partial X_l \partial X_m\partial X_n\partial X_p\partial X_q\partial X_r} H_lH_mH_nH_pH_qH_r\\ +\dfrac 1 {5!} \sum_{m,n,p,q,r} \dfrac{\partial^5 Q_{23}^i}{\partial X_m\partial X_n\partial X_p\partial X_q\partial X_r} H_m H_nH_pH_qH_r =0
\end{array}$$

$$\begin{array}{l} Q_{38}^i=\dfrac 1 {7!} \sum_{l,m,n,p,q,r,s} \dfrac{\partial^7 Q_{22}^i}{ \partial X_l \partial X_m\partial X_n\partial X_p\partial X_q\partial X_r\partial X_s} H_lH_mH_nH_pH_qH_rH_s\\ +\dfrac 1 {6!} \sum_{m,n,p,q,r,s} \dfrac{\partial^6 Q_{23}^i}{\partial X_m\partial X_n\partial X_p\partial X_q\partial X_r\partial X_s} H_m H_nH_pH_qH_rH_s =0
\end{array}$$

\noindent Using (\ref{eq8}) and (\ref{eq9}),

$$ Q_{39}^i=\dfrac 1 {7!} \sum_{m,n,p,q,r,s,t} \dfrac{\partial^7 Q_{23}^i}{\partial X_m\partial X_n\partial X_p\partial X_q\partial X_r\partial X_s\partial X_t} H_m H_nH_pH_qH_rH_sH_t =0
$$

\noindent Using (\ref{eq9}),

$$ Q_{3,10}^i=\dfrac 1 {8!} \sum_{m,n,p,q,r,s,t,u} \dfrac{\partial^8 Q_{23}^i}{\partial X_m\partial X_n\partial X_p\partial X_q\partial X_r\partial X_s\partial X_t\partial X_u} H_m H_nH_pH_qH_rH_sH_tH_u =0
$$

\noindent And using (\ref{eq9}) and (\ref{eq10}),

$$ Q_{3,11}^i=\dfrac 1 {9!} \sum_{m,n,p,q,r,s,t,u,v} \dfrac{\partial^9 Q_{23}^i}{\partial X_m\partial X_n\partial X_p\partial X_q\partial X_r\partial X_s\partial X_t\partial X_u\partial X_v} H_m H_nH_pH_qH_rH_sH_tH_uH_v =0
$$

\noindent We have then obtained

$$\begin{array}{lll} P_3^i &=& Q_{31}^i+Q_{32}^i\\ &=& 6L_i \sum_{j,k=1}^d a_{ij}a_{ik} L_j^3 L_k^3+\sum_{j,k,l=1}^d 6 a_{ij} a_{ik} a_{il} L_j^3L_k^3L_l^3. \end{array} $$

\noindent In particular, $\deg P_3^i \leq 9$. The expression
obtained for $P_3^i$ gives

$$P_4^i= \sum_{j=1}^{10} Q_{4j}^i \text{\ with \ } \deg Q_{4j}^i=2j+7.$$

In order to determine $P_4^i$ we need to compute the derivatives of $Q_{31}^i$ and $Q_{32}^i$. We compute first the derivatives of $Q_{31}^i$.

$$\dfrac{\partial Q_{31}^i}{\partial X_l}  = \sum_{j,k=1}^d(6 a_{ij}a_{ik}a_{il}  L_j^3 L_k^3+18 a_{ij}a_{ik} a_{jl} L_i  L_j^2 L_k^3+18 a_{ij} a_{ik}a_{kl} L_i  L_j^3 L_k^2)$$

$$\begin{array}{lll} \dfrac{\partial^2 Q_{31}^i}{\partial X_l\partial X_m}  &=& \sum_{j,k=1}^d(18 a_{ij}a_{ik}a_{il} a_{jm} L_j^2 L_k^3+18 a_{ij}a_{ik}a_{il} a_{km} L_j^3 L_k^2+18 a_{ij}a_{ik} a_{jl}a_{im} L_j^2 L_k^3\\ && +36 a_{ij}a_{ik} a_{jl}a_{jm} L_i  L_j L_k^3+54 a_{ij}a_{ik} a_{jl}a_{km} L_i  L_j^2 L_k^2+18 a_{ij} a_{ik}a_{kl} a_{im}  L_j^3 L_k^2 \\ && +54 a_{ij} a_{ik}a_{kl} a_{jm} L_i  L_j^2 L_k^2+36 a_{ij} a_{ik}a_{kl} a_{km}L_i  L_j^3 L_k).\end{array} $$

\noindent The derivative $\partial^3 Q_{31}^i/\partial X_l\partial X_m\partial X_n$ is a sum of terms of one of the following forms

$$\begin{array}{ll} 36 \sum_{j,k=1}^da_{ij}a_{ik}a_{il} a_{jm} a_{jn} L_j L_k^3, & 54 \sum_{j,k=1}^da_{ij}a_{ik}a_{il} a_{jm} a_{kn} L_j^2 L_k^2 \\ 36 \sum_{j,k=1}^da_{ij}a_{ik} a_{jl}a_{jm} a_{jn} L_i  L_k^3, & 108 \sum_{j,k=1}^da_{ij}a_{ik} a_{jl}a_{jm} a_{kn} L_i  L_j L_k^2, \end{array}$$

\noindent up to a swift of $j$ with $k$ and a permutation of the set $\{l,m,n\}$. The derivative $\partial^4 Q_{31}^i/\partial X_l\partial X_m\partial X_n\partial X_p$ is a sum of terms of one of the following forms

$$\begin{array}{ll} 36 \sum_{j,k=1}^d a_{ij}a_{ik}a_{il} a_{jm} a_{jn} a_{jp} L_k^3, &  108 \sum_{j,k=1}^d a_{ij}a_{ik}a_{il} a_{jm} a_{jn} a_{kp} L_j L_k^2\\ 108 \sum_{j,k=1}^da_{ij}a_{ik} a_{jl}a_{jm} a_{jn} a_{kp} L_i  L_k^2, & 216 \sum_{j,k=1}^da_{ij}a_{ik} a_{jl}a_{jm} a_{kn}a_{kp} L_i  L_j L_k, \end{array}$$

\noindent up to a swift of $j$ with $k$ and a permutation of the set $\{l,m,n,p\}$. The derivative $\partial^5 Q_{31}^i/\partial X_l\partial X_m\partial X_n\partial X_p \partial X_q$ is a sum of terms of one of the following forms

$$\begin{array}{ll} 108 \sum_{j,k=1}^d a_{ij}a_{ik}a_{il} a_{jm} a_{jn} a_{jp} a_{kq} L_k^2, &  216 \sum_{j,k=1}^d a_{ij}a_{ik}a_{il} a_{jm} a_{jn} a_{kp}a_{kq} L_j L_k\\ 216 \sum_{j,k=1}^d a_{ij}a_{ik} a_{jl}a_{jm} a_{jn} a_{kp} a_{kq} L_i  L_k, &  \end{array}$$

\noindent up to a swift of $j$ with $k$ and a permutation of the set $\{l,m,n,p,q\}$. The derivative $\partial^6 Q_{31}^i/\partial X_l\partial X_m\partial X_n\partial X_p \partial X_q\partial X_r$ is a sum of terms of one of the following forms

$$ 216 \sum_{j,k=1}^d a_{ij}a_{ik}a_{il} a_{jm} a_{jn} a_{jp} a_{kq}a_{kr} L_k, \quad  216 \sum_{j,k=1}^d a_{ij}a_{ik} a_{jl}a_{jm} a_{jn} a_{kp} a_{kq} a_{kr} L_i,$$

\noindent up to a swift of $j$ with $k$ and a permutation of the set $\{l,m,n,p,q,r\}$. The derivative $\partial^7 Q_{31}^i/\partial X_l\partial X_m\partial X_n\partial X_p \partial X_q\partial X_r\partial X_s$ is a sum of terms of the following form

$$ 216 \sum_{j,k=1}^d a_{ij}a_{ik}a_{il} a_{jm} a_{jn} a_{jp} a_{kq}a_{kr}a_{ks},$$

\noindent up to a swift of $j$ with $k$ and a permutation of the set $\{l,m,n,p,q,r,s\}$.

We compute now the derivatives of $Q_{32}^i$.

$$\dfrac{\partial Q_{32}^i}{\partial X_m} = 18 \sum_{j,k,l}(  a_{ij} a_{ik} a_{il}a_{jm} L_j^2L_k^3L_l^3+  a_{ij} a_{ik} a_{il} a_{km} L_j^3L_k^2L_l^3+a_{ij} a_{ik} a_{il} a_{lm} L_j^3L_k^3L_l^2)$$

\noindent The derivative $\partial^2 Q_{32}^i/\partial X_m\partial X_n$ is a sum of terms of one of the following forms

$$36 \sum_{j,k,l}a_{ij} a_{ik} a_{il}a_{jm} a_{jn} L_jL_k^3L_l^3, \quad 54 \sum_{j,k,l} a_{ij} a_{ik} a_{il}a_{jm} a_{kn} L_j^2L_k^2L_l^3$$

\noindent up to a permutation of the set $\{j,k,l\}$ and a permutation of the set $\{m,n\}$. The derivative $\partial^3 Q_{32}^i/\partial X_m\partial X_n\partial X_p$ is a sum of terms of one of the following forms

$$\begin{array}{ll} 36 \sum_{j,k,l} a_{ij} a_{ik} a_{il}a_{jm} a_{jn} a_{jp} L_k^3L_l^3, & 108 \sum_{j,k,l} a_{ij} a_{ik} a_{il}a_{jm} a_{jn} a_{kp} L_jL_k^2L_l^3 \\ 162 \sum_{j,k,l} a_{ij} a_{ik} a_{il}a_{jm} a_{kn}a_{lp}  L_j^2L_k^2L_l^2, &  \end{array}$$

\noindent up to a permutation of the set $\{j,k,l\}$ and a permutation of the set $\{m,n,p\}$. The derivative $\partial^4 Q_{32}^i/\partial X_m\partial X_n\partial X_p\partial X_q$ is a sum of terms of one of the following forms

$$\begin{array}{ll} 108 \sum_{j,k,l} a_{ij} a_{ik} a_{il}a_{jm} a_{jn} a_{jp}a_{kq}  L_k^2L_l^3, & 216 \sum_{j,k,l} a_{ij} a_{ik} a_{il}a_{jm} a_{jn} a_{kp} a_{kq} L_jL_k L_l^3 \\ 324 \sum_{j,k,l} a_{ij} a_{ik} a_{il}a_{jm} a_{jn} a_{kp}a_{lq}  L_jL_k^2L_l^2, &  \end{array}$$

\noindent up to a permutation of the set $\{j,k,l\}$ and a permutation of the set $\{m,n,p,q\}$. The derivative $\partial^5 Q_{32}^i/\partial X_m\partial X_n\partial X_p\partial X_q\partial X_r$ is a sum of terms of one of the following forms

$$\begin{array}{ll} 216 \sum_{j,k,l} a_{ij} a_{ik} a_{il}a_{jm} a_{jn} a_{jp}a_{kq} a_{kr}  L_k L_l^3, & 324 \sum_{j,k,l} a_{ij} a_{ik} a_{il}a_{jm} a_{jn} a_{jp}a_{kq}a_{lq}   L_k^2L_l^2 \\  648 \sum_{j,k,l} a_{ij} a_{ik} a_{il}a_{jm} a_{jn} a_{kp} a_{kq}a_{lr}  L_jL_k L_l^2, &  \end{array}$$

\noindent up to a permutation of the set $\{j,k,l\}$ and a permutation of the set $\{m,n,p,q,r\}$. The derivative $\partial^6 Q_{32}^i/\partial X_m\partial X_n\partial X_p\partial X_q\partial X_r\partial X_s$ is a sum of terms of one of the following forms

$$\begin{array}{ll} 216 \sum_{j,k,l} a_{ij} a_{ik} a_{il}a_{jm} a_{jn} a_{jp}a_{kq} a_{kr}a_{ks}  L_l^3, & 648 \sum_{j,k,l} a_{ij} a_{ik} a_{il}a_{jm} a_{jn} a_{jp}a_{kq}a_{lq} a_{ks}  L_kL_l^2 \\  1296 \sum_{j,k,l} a_{ij} a_{ik} a_{il}a_{jm} a_{jn} a_{kp} a_{kq}a_{lr} a_{ls} L_jL_k L_l, &  \end{array}$$

\noindent up to a permutation of the set $\{j,k,l\}$ and a permutation of the set $\{m,n,p,q,r,s\}$. The derivative $\partial^7 Q_{32}^i/\partial X_m\partial X_n\partial X_p\partial X_q\partial X_r\partial X_s\partial X_t$ is a sum of terms of one of the following forms

$$ 648 \sum_{j,k,l} a_{ij} a_{ik} a_{il}a_{jm} a_{jn} a_{jp}a_{kq} a_{kr}a_{ks} a_{lt} L_l^2, \quad 1296 \sum_{j,k,l} a_{ij} a_{ik} a_{il}a_{jm} a_{jn} a_{jp}a_{kq}a_{lq} a_{ks}a_{lt}  L_kL_l $$

\noindent up to a permutation of the set $\{j,k,l\}$ and a permutation of the set $\{m,n,p,q,r,s,t\}$. The derivative $\partial^8 Q_{32}^i/\partial X_m\partial X_n\partial X_p\partial X_q\partial X_r\partial X_s\partial X_t\partial X_u$ is a sum of terms of the form

$$ 1296 \sum_{j,k,l} a_{ij} a_{ik} a_{il}a_{jm} a_{jn} a_{jp}a_{kq} a_{kr}a_{ks} a_{lt}a_{lu} L_l,$$

\noindent up to a permutation of the set $\{j,k,l\}$ and a permutation of the set $\{m,n,p,q,r,s,t,u\}$. The derivative $\partial^9 Q_{32}^i/\partial X_m\partial X_n\partial X_p\partial X_q\partial X_r\partial X_s\partial X_t\partial X_u\partial X_v$ is a sum of terms of the form

$$ 1296 \sum_{j,k,l} a_{ij} a_{ik} a_{il}a_{jm} a_{jn} a_{jp}a_{kq} a_{kr}a_{ks} a_{lt}a_{lu}a_{lv},$$

\noindent up to a permutation of the set $\{j,k,l\}$ and a permutation of the set $\{m,n,p,q,r,s,t,u,v\}$.

We compute now the homogeneous summands of $P_4^i$.

Using (\ref{eq8}),

$$Q_{41}^i = \sum_{l=1}^d \dfrac{\partial Q_{31}^i}{\partial X_l}H_l = \sum_{j,k,l=1}^d 6 a_{ij}a_{ik}a_{il}  L_j^3 L_k^3 L_l^3.$$

\noindent Using (\ref{eq8}), (\ref{eq9}) and (\ref{eq10}), we obtain

$$Q_{4j}^i=0, \quad \forall j=2,\dots,10.$$

\noindent We have then

$$P_4^i= Q_{41}^i = \sum_{j,k,l=1}^d 6 a_{ij}a_{ik}a_{il}  L_j^3 L_k^3 L_l^3,$$

\noindent hence $P_4^i$ is a homogeneous polynomial of degree 9.
The expression obtained for $P_4^i$ gives

$$P_5^i= \sum_{j=1}^9 Q_{5j}^i \text{\ with \ } \deg Q_{5j}^i=2j+9.$$

\noindent In order to determine $P_5^i$ we compute the derivatives
of $Q_{41}^i$.

$$\dfrac{\partial Q_{41}^i}{\partial X_m}= 18 \sum_{j,k,l=1}^d(a_{ij}a_{ik}a_{il} a_{jm} L_j^2 L_k^3 L_l^3+a_{ij}a_{ik}a_{il} a_{km} L_j^3 L_k^2 L_l^3+a_{ij}a_{ik}a_{il}a_{lm}  L_j^3 L_k^3 L_l^2).$$

\noindent The derivative $\partial^2 Q_{41}^i/\partial X_m \partial X_n$ is a sum of terms of one of the following forms

 $$36 \sum_{j,k,l=1}^d  a_{ij}a_{ik}a_{il} a_{jm}a_{jn} L_j L_k^3 L_l^3, \quad 54 \sum_{j,k,l=1}^d a_{ij}a_{ik}a_{il} a_{jm}a_{kn} L_j^2 L_k^2 L_l^3,$$

\noindent up to a permutation of the set $\{j,k,l\}$ and a permutation of the set $\{m,n\}$. The derivative $\partial^3 Q_{41}^i/\partial X_m \partial X_n\partial X_p$ is a sum of terms of one of the following forms

 $$\begin{array}{ll} 36 \sum_{j,k,l=1}^d  a_{ij}a_{ik}a_{il} a_{jm}a_{jn} a_{jp} L_k^3 L_l^3, & 108 \sum_{j,k,l=1}^d  a_{ij}a_{ik}a_{il} a_{jm}a_{jn} a_{kp} L_j L_k^2 L_l^3,\\ 162 \sum_{j,k,l=1}^d a_{ij}a_{ik}a_{il} a_{jm}a_{kn}a_{lp} L_j^2 L_k^2 L_l^2, & \end{array} $$

\noindent up to a permutation of the set $\{j,k,l\}$ and a permutation of the set $\{m,n,p\}$. The derivative $\partial^4 Q_{41}^i/\partial X_m \partial X_n\partial X_p\partial X_q$ is a sum of terms of one of the following forms

 $$\begin{array}{ll} 108 \sum_{j,k,l=1}^d  a_{ij}a_{ik}a_{il} a_{jm}a_{jn} a_{jp} a_{kq} L_k^2 L_l^3, & 216 \sum_{j,k,l=1}^d  a_{ij}a_{ik}a_{il} a_{jm}a_{jn} a_{kp} a_{kq} L_j L_k L_l^3,\\ 324 \sum_{j,k,l=1}^d a_{ij}a_{ik}a_{il} a_{jm}a_{kn}a_{lp}a_{jq} L_j L_k^2 L_l^2, & \end{array} $$

\noindent up to a permutation of the set $\{j,k,l\}$ and a permutation of the set $\{m,n,p,q\}$. The derivative $\partial^5 Q_{41}^i/\partial X_m \partial X_n\partial X_p\partial X_q\partial X_r$ is a sum of terms of one of the following forms

 $$\begin{array}{ll} 216 \sum_{j,k,l=1}^d  a_{ij}a_{ik}a_{il} a_{jm}a_{jn} a_{jp} a_{kq}a_{kr}  L_k L_l^3, & 324 \sum_{j,k,l=1}^d  a_{ij}a_{ik}a_{il} a_{jm}a_{jn} a_{jp} a_{kq} a_{lr} L_k^2 L_l^2, \\ 648 \sum_{j,k,l=1}^d  a_{ij}a_{ik}a_{il} a_{jm}a_{jn} a_{kp} a_{kq}a_{lr} L_j L_k L_l^2, & \end{array} $$

\noindent up to a permutation of the set $\{j,k,l\}$ and a permutation of the set $\{m,n,p,q,r\}$. The derivative $\partial^6 Q_{41}^i/\partial X_m \partial X_n\partial X_p\partial X_q\partial X_r\partial X_s$ is a sum of terms of one of the following forms

 $$\begin{array}{ll} 648 \sum_{j,k,l=1}^d  a_{ij}a_{ik}a_{il} a_{jm}a_{jn} a_{jp} a_{kq}a_{kr} a_{ls} L_k L_l^2, & 216 \sum_{j,k,l=1}^d  a_{ij}a_{ik}a_{il} a_{jm}a_{jn} a_{jp} a_{kq}a_{kr}a_{ks} L_l^3, \\ 1296 \sum_{j,k,l=1}^d  a_{ij}a_{ik}a_{il} a_{jm}a_{jn} a_{kp} a_{kq}a_{lr} a_{ls} L_j L_k L_l, & \end{array} $$

\noindent up to a permutation of the set $\{j,k,l\}$ and a permutation of the set $\{m,n,p,q,r,s\}$. The derivative $\partial^7 Q_{41}^i/\partial X_m \partial X_n\partial X_p\partial X_q\partial X_r\partial X_s\partial X_t$ is a sum of terms of one of the following forms

 $$648 \sum_{j,k,l=1}^d  a_{ij}a_{ik}a_{il} a_{jm}a_{jn} a_{jp} a_{kq}a_{kr} a_{ls}a_{kt} L_l^2, \, 1296 \sum_{j,k,l=1}^d  a_{ij}a_{ik}a_{il} a_{jm}a_{jn} a_{jp} a_{kq}a_{kr} a_{ls}a_{lt} L_k L_l, $$

\noindent up to a permutation of the set $\{j,k,l\}$ and a permutation of the set $\{m,n,p,q,r,s,t\}$. The derivative $\partial^8 Q_{41}^i/\partial X_m \partial X_n\partial X_p\partial X_q\partial X_r\partial X_s\partial X_t\partial X_u$ is a sum of terms of the following form

 $$1296 \sum_{j,k,l=1}^d  a_{ij}a_{ik}a_{il} a_{jm}a_{jn} a_{jp} a_{kq}a_{kr} a_{ls}a_{kt}a_{lu} L_l, $$

\noindent up to a permutation of the set $\{j,k,l\}$ and a permutation of the set $\{m,n,p,q,r,s,t,u\}$. The derivative $\partial^9 Q_{41}^i/\partial X_m \partial X_n\partial X_p\partial X_q\partial X_r\partial X_s\partial X_t\partial X_u\partial X_v$ is a sum of terms of the following form

 $$1296 \sum_{j,k,l=1}^d  a_{ij}a_{ik}a_{il} a_{jm}a_{jn} a_{jp} a_{kq}a_{kr} a_{ls}a_{kt}a_{lu} a_{lv}, $$

\noindent up to a permutation of the set $\{j,k,l\}$ and a permutation of the set $\{m,n,p,q,r,s,t,u,v\}$. We have

$$Q_{5j}^i= \dfrac 1 {j!} \sum_{n_1,\dots,n_j=1}^d \dfrac{\partial^j Q_{41}^i}{\partial X_{n_1} \dots \partial X_{n_j}} H_{n_1}\dots H_{n_j}$$

\noindent and, using again (\ref{eq8}), (\ref{eq9}), (\ref{eq10}), we obtain

$$Q_{5j}^i =0, \text{\ for all \ }j,$$

\noindent hence $P_5^i=0$, for all $i=1\dots,d$, as wanted. \hfill $\Box$

\vspace{0.5cm}
The following example shows that the result in Proposition \ref{prop} is optimal, in the sense that, under the hypothesis $(JH)^3=0$, we may not expect that the first zero term in the sequence $(P_j^i)$ comes before the fifth for all $i=1,\dots,d$.

\begin{example} We consider the polynomial map $F$ as in Theorem \ref{theo2} with $d=5$ and

$$\left\{ \begin{array}{lll} L_1(X_1,X_2,X_3,X_4,X_5)&=& a_2X_2+a_3X_3+a_4X_4+a_5X_5 \\ L_2(X_1,X_2,X_3,X_4,X_5)&=& b_3X_3+b_4X_4+b_5X_5\\ L_3(X_1,X_2,X_3,X_4,X_5&=& 0 \\  L_4(X_1,X_2,X_3,X_4,X_5)&=& 0\\  L_5(X_1,X_2,X_3,X_4,X_5)&=& 0  \end{array} \right.$$

\noindent with $a_2,a_3,a_4,a_5,b_3,b_4,b_5$ complex parameters.
We may check that the Jacobian matrix $JH$ is nilpotent of degree
3. By computing the sequences of polynomials $P_1^i$, we obtain

$$\begin{array}{lll} P_2^1 &=& 3 a_2 L_1^2 L_2^3 +3 a_2^2 L_1 L_2^6 + a_2^3 L_2^9 \\[8pt] P_3^1 &=& 6 a_2^2 L_1L_2^6+6a_2^3L_2^9\\[8pt]
P_4^1&=& 6a_2^3 L_2^9 \\ [8pt] P_5^1&=&0 \end{array} $$

\noindent Clearly, $P_2^2=0$ and $P_1^3=P_1^4=P_1^5=0$. Then,
choosing such an $F$ with $a_2 \neq 0$ and $L_2$ not identically
zero, we obtain an example of a polynomial map $F$ such that the
Jacobian matrix $JH$ is nilpotent of degree 3 and at least one
$P_4^i$ is not zero. Hence the inverse map $F^{-1}$ has degree
equal to $9$. The map $G=F^{-1}$ is given by

$$\left\{ \begin{array}{lll} G_1(Y) &=& Y_1-L_1(Y)^3+3 a_2 L_1(Y)^2 L_2(Y)^3 -3 a_2^2 L_1(Y) L_2(Y)^6 + a_2^3 L_2(Y)^9\\
G_2(Y) &=& Y_2-L_2(Y)^3\\ G_3(Y) &=& Y_3
\\ G_4(Y) &=& Y_4\\
G_5(Y) &=& Y_5 \end{array} \right. $$

\noindent where $Y:=(Y_1,Y_2,Y_3,Y_4,Y_5)$.
\end{example}

\vspace{0.3cm}
The result obtained in Theorem \ref{theo2} leads us
to formulate the following conjecture.

\begin{conjecture}\label{conj}
{\it Let $F:K^d \rightarrow K^d$ be a polynomial map of the form

$$\left\{ \begin{array}{lll} F_1(X_1,\dots,X_d)&=& X_1+H_1(X_1,\dots,X_d) \\ F_2(X_1,\dots,X_d)&=& X_2+H_2(X_1,\dots,X_d)\\ & \vdots & \\  F_d(X_1,\dots,X_d)&=& X_d+H_d(X_1,\dots,X_d) , \end{array} \right.$$

\noindent where $H_i(X_1,\dots,X_d)=L_i(X_1,\dots,X_d)^3$ and $L_i(X_1,\dots,X_d)=a_{i1}X_1+\dots +a_{id} X_d$, $1\leq i \leq d$. We consider the Jacobian matrix $JH$ of $H:=(H_1,\dots,H_d)$. For each $i=1,\dots,d$, we consider the polynomial sequence $(P_j^i)$ defined in the following way

$$\begin{array}{lll} P_0^i(X_1,\dots,X_d) &= & X_i, \\
P_1^i(X_1,\dots,X_d) &=& H_i , \end{array}$$

\noindent and, assuming $P_{j-1}^i$ is defined,

$$P_j^i(X_1,\dots,X_d) = P_{j-1}^i(F_1,\dots,F_d)-P_{j-1}^i(X_1,\dots,X_d).$$

Let $g$ be an integer, $1\leq g \leq d$. If $(JH)^g=0$, then $P_{(3^{g-1}+1)/2}^i=0$, for all $i=1,\dots,d$, $F$ is invertible and the inverse of $F$ has maximal degree at most equal to $3^{g-1}$.}
\end{conjecture}

\vspace{0.3cm}

\begin{remark} By Theorem \ref{theo}, the condition $P_{(3^{g-1}+1)/2}^i=0$, for all $i=1,\dots,d$, implies that $F$ is invertible.
Since, for $F$ as in Conjecture \ref{conj}, $\det(JF)=1$ implies
$(JH)^d=0$, and taking into account Proposition \ref{red},
Conjecture \ref{conj} implies the Jacobian conjecture for complex polynomial maps.
\end{remark}

\begin{remark} Conjecture \ref{conj} is true for

\begin{enumerate}[1)]
\item $g=1$: since $H_i$ is homogeneous, $JH=0$ implies $H=0$, hence $P_1^i=0, \, \forall i=1,\dots,d$ and $F=Id$.
\item $g=2$: if $(JH)^2=0$, we have $\sum_{j=1}^d \dfrac{\partial H_i}{\partial X_j}\dfrac{\partial H_j}{\partial X_k}=0, \, \forall \, i, k =1,\dots,d.$
Multiplying by $X_k$ and summing up from $k=1$ to $d$, we obtain
$\sum_{j=1}^d \dfrac{\partial H_i}{\partial X_j}H_j=0, \, \forall
\, i$. Then $F$ is a quasi-translation, $P_2^i=0$, for all
$i=1,\dots,d$, and $F^{-1}=X-H$ (see \cite{Alg} Proposition 9).
\item $g=3$: this is Proposition \ref{prop}.
\end{enumerate}
\end{remark}

\end{document}